\documentclass[reqno]{amsart}
\usepackage[utf8x]{inputenc}
\usepackage{combelow}
\usepackage{amsmath}
\usepackage{amssymb}
 \usepackage{paralist}
 \usepackage{graphics} %% add this and next lines if pictures should be in esp format
 \usepackage{epsfig} %For pictures: screened artwork should be set up with an 85 or 100 line screen
\usepackage{graphicx}  \usepackage{epstopdf}%This is to transfer .eps figure to .pdf figure; please compile your paper using PDFLeTex or PDFTeXify.
 \usepackage[colorlinks=true]{hyperref}
\usepackage{graphicx}
\usepackage{epsf}
\usepackage{float}
\usepackage[dvipsnames]{xcolor}
\usepackage{booktabs}
\usepackage{color}
\usepackage{subfigure}
\usepackage{multicol}
\usepackage{enumerate}

\newcommand{\be}{\begin{equation}}
\newcommand{\ee}{\end{equation}}
\newcommand{\bea}{\begin{eqnarray}}
\newcommand{\eea}{\end{eqnarray}}
\newcommand{\beas}{\begin{eqnarray*}}
\newcommand{\eeas}{\end{eqnarray*}}
\newcommand{\bl}{\begin{law}}
\newcommand{\el}{\end{law}}
\newcommand{\bthm}{\begin{thm}}
\newcommand{\ethm}{\end{thm}}

\newcommand{\ud}{\,\mathrm{d}}

\newcommand{\G}{\mathrm{G}}

\newcommand{\I}{\mathrm{i}}

\newcommand{\mR}{\mathbb R}

\newcommand{\barint}{\kern4pt \raise3.4pt\hbox{\vrule height.6pt
    width7pt} \kern-11pt \int}

\newcommand{\erf}{\text{erf}}

\newcommand{\Cc}{\mathrm{C}}
\newcommand{\Ss}{\mathrm{S}}

\newcommand{\ds}{\displaystyle}

\usepackage{multirow}
\reversemarginpar
\newcounter{mnote}
\hypersetup{urlcolor=blue, citecolor=red}

\title[Fresnel Integral Computation Techniques]
      {Fresnel Integral Computation Techniques}

\author[A. Ionu\cb{t}, J. C. Hateley]{Alexandru Ionu\cb{t}, James C. Hateley}
\email{alexandru.ionut172@gmail.com, hateleyjc@gmail.com}

\begin{document}

% Enter the first author's name and address:
%{\footnotesize
% please put the address of the first author
% \centerline{Department of Mathematics}
%   \centerline{Vanderbilt University}
%   \centerline{Nashville, TN 37027}
%} % Do not forget to end the {\footnotesize by the sign }

%\bigskip

% The name of the associate editor will be entered by an editorial staff
% "Communicated by the associate editor name" is not needed for special issue.
% \centerline{(Communicated by the associate editor name)}

%The abstract of your paper
\begin{abstract}
This work is an extension of previous work by Alazah et al. [M. Alazah, S. N. Chandler-Wilde, and S. La Porte,  Numerische Mathematik, 128(4):635–661, 2014]. We split the computation of the Fresnel Integrals into 3 cases: a truncated Taylor series, modified trapezoid rule and an asymptotic expansion for small, medium and large arguments respectively. These special functions can be computed accurately and efficiently up to an arbitrary precision.  Error estimates are provided and we give a systematic method in choosing the various parameters for a desired precision. We illustrate this method and verify numerically using double precision.
\end{abstract}

\maketitle
%%%%%%%%%%%%%%%%%%%%
%---CONTENT----------------------------------------------------------------------------------------
%---INTRODUCTION-----------------------------------------------------------------------------------
\section{Introduction}\label{sec:1}
The Fresnel integrals and their simultaneous parametric plot, the clothoid, have numerous applications including; but not limited to, optics and electromagnetic theory \cite{goloborodko2016computer,sandoval2018approximation,sangeorzan2008some,shimobaba2008numerical,zhang2016layered}, robotics \cite{bevilacqua2016path,fleury1995primitives,montes2008real,nagy2001trajectory,reuter1998mobile}, civil engineering \cite{butz2004realistic,kim2010comparison,wang2001approximation} and biology \cite{PMID:31998835}. There have been numerous works over the past 70 years computing and numerically approximating Fresnel integrals. Boersma established approximations using the Lanczos tau-method~\cite{Boersma1960ComputationOF} and Cody computed rational Chebyshev approximations using the Remes algorithm \cite{cody1968chebyshev}. Another approach includes a spreadsheet computation by Mielenz~\cite{mielenz1997computation}; which is based on successive improvements of known relational approximations.  Mielenz also gives an improvement of his work~\cite{mielenz2000computation}, where the accuracy is less then $1.e$-$9$. More recently, Alazah, Chandler-Wilde and LaPorte propose a method to compute these integrals via a modified trapezoid rule~\cite{alazah2014computing}.

Alazah et al. remark after some experimentation that a truncation of the Taylor series is more efficient and accurate than their new method for a small argument~\cite{alazah2014computing}. We build on this philosophy and add to it by endorsing asymptotic expansion for large arguments (see~\cite{mielenz1997computation} for an earlier example). We propose a technique combining truncated Taylor series, the modified trapezium rule and asymptotic expansion and test its validity for practical computation of the Fresnel integrals all over the real line. Cut-off points are found analytically yielding a technique valid for arbitrary precision. 

The paper is organized as follows: Section~\ref{sec:2} gives a brief introduction to Fresnel integrals. In Section~\ref{sec:3} the approximations with appropriate bounds  are introduced. Section~\ref{sec:4} contains a systematic method for choosing the parameters to construct an approximation using results from Section~\ref{sec:3}.  Section~\ref{sec:5} contains some benchmarking and comparative results.  Lastly, conclusions are provided in Section~\ref{sec:6}. 
%---PRELIMINARY------------------------------------------------------------------------------------
\section{Fresnel Integrals}\label{sec:2}
We define the Fresnel cosine and sine integrals, $\Cc(x)$ and $\Ss(x)$ respectively, as follows:
\begin{equation}
\Cc(x)=\int_{0}^{x} \cos\Big(\frac{\pi t^2}{2} \Big)\ud t, \qquad
\Ss(x)=\int_{0}^{x} \sin\Big(\frac{\pi t^2}{2} \Big)\ud t .
\end{equation}
This is consistent with~\cite{alazah2014computing}. With this pair of integrals, $\Cc(x)$ and $\Ss(x)$ can be written as a complex exponential integral:
\begin{equation}
\G(x) = \Cc(x) +\I \Ss(x) = \int_{0}^{x} \exp\Big(\frac{\I\pi t^2}{2}\Big)\ud t.
\end{equation}
The Argand diagram of $\G(x)$ for real $x$ plots the clothoid. Also, both $\Cc(x)$ and $\Ss(x)$ are odd functions of $x$ and can be extended to analytic functions over the complex plane. %This is easy to see as one can write
%\begin{align*}
%\Cc(x) =& \ds\frac{\sqrt{\pi}}{4}\left(\sqrt{\I}\erf(\sqrt{\I}x) + \sqrt{-\I}\erf(\sqrt{-\I}x) \right) \\
%\Ss(x) =& \ds\frac{\sqrt{\pi}}{4}\left(\sqrt{-\I}\erf(\sqrt{\I}x) + \sqrt{\I}\erf(\sqrt{-\I}x) \right)
%\end{align*}
$\G(x)$ can also be expressed as
\begin{equation}
\G(x) = \frac{1+i}{2}\erf\Big(\frac{\sqrt{\pi}}{2}(1-\I)x\Big)
\end{equation}
with $\erf(x) = \frac{2}{\pi}\int_0^x \exp(-t^2)\ud t$ being the error function.
%---APPROXIMATION----------------------------------------------------------------------------------
\section{Approximations and Error bounds}\label{sec:3}
In this section we list the approximations to be used for $\mathrm{G}(x)$ along with their corresponding error bounds.
We start with the truncated Taylor series $\mathrm{T}_N(x)$;
\begin{equation}\label{eq:TN}
\mathrm{T}_N(x) = \sum_{k=0}^{N} \frac{(\I\pi)^k x^{2k+1}}{2^k(2k+1)k!}
\end{equation}
$T_N(x)$ is an alternating series in both real and imaginary parts bounding these by the next term, we have the estimate
\begin{equation}\label{eq:ETN}
|\G(x) - \mathrm{T}_N(x)| \leq E_{\mathrm{T}_N}(x) = \frac{\pi^{N+1} x^{2N+3}}{2^{N+1}(2N+3)(N+1)!} + \frac{\pi^{N+2} x^{2N+5}}{2^{N+2}(2N+5)(N+2)!}
\end{equation}
Modified trapezoid sum from \cite{alazah2014computing}, which is
%\begin{multline}\label{eq:GN}
\begin{align}\label{eq:GN}
%\frac{1+\I}{2}-\frac{1+\I}{\exp\left((1-\I)\pi A_N x\right)+1}
\mathrm{G}_N(x) = \frac{1+\I}{2}- &\frac{1+\I}{\exp\left((1-\I)\pi A_N x\right)+1} \\ 
&-\frac{2\I x\exp\big(\I\pi x^2/2\big)}{\pi A_N}\sum_{k=1}^{N} \frac{\exp\left(-\pi(k-1/2)^2A_N^{-2}\right)}{x^2+\I 2(k-1/2)^2 A_N^{-2}},
%\end{multline}
\end{align}
where $A_N$ is defined in \eqref{eq:GN:consants}. A global bound for eq. \eqref{eq:GN}, see eq. (67) from \cite{alazah2014computing}, is
\begin{equation}\label{eq:EGN}
|\mathrm{G}(x) - \mathrm{G}_N(x)| \leq E_{\mathrm{G}_N} = \ds\frac{2\sqrt{2}c_N\exp(-\pi N)}{2N+1}.
\end{equation}
Where $c_N$ is given as eq. (47) in \cite{alazah2014computing}; which is
\begin{equation*}
c_N = \ds\frac{20 \sqrt{2} \exp(-\pi/2)}{9\pi\big(1-\exp(-2\pi A_N^2)\big)}\big(1+2\sqrt{\pi}\exp(-\beta\pi A_N^2)\big) + \ds\frac{(2\pi+1)\exp(-\pi/2)}{2\sqrt{2}\pi^{3/2}}
\end{equation*}
and $\beta$, $A_N$ are the constants
\begin{equation}\label{eq:GN:consants}
\beta = 1 - \ds\frac{1}{\sqrt{2}} - \frac{2\sqrt{2}+ 1}{16}, \qquad A_N = \sqrt{N + 1/2}.
\end{equation}
The asymptotic expansion can be derived directly as follows, first write
\begin{equation*}
\mathrm{G}(x) = \frac{1+\I}{2} - \int_x^\infty \exp\Big(\frac{\I\pi t^2}{2}\Big)\ud t.
\end{equation*}
Applying repeated applications of integration by parts gives
\begin{equation}\label{eq:QN:ibt} 
\mathrm{G}(x) = \mathrm{Q}_N(x) + \int_x^\infty \frac{(2N+1)!!(-\I)^{N+1}}{\pi^{N+1} t^{2k+2}}\exp\Big(\frac{\I\pi t^2}{2}\Big) \ud t
\end{equation}
with
\begin{equation}\label{eq:QN}
\mathrm{Q}_N(x) = \frac{1+\I}{2} + \exp\Big(\frac{\I\pi x^2}{2}\Big) \sum_{k=0}^{N} \frac{(2k-1)!!(-\I)^{k+1}}{\pi^{k+1} x^{2k+1}} 
\end{equation}
For the error bound, estimating the  integral in \eqref{eq:QN:ibt} gives,
\begin{equation}\label{eq:EQN}
%|G(x) - Q_N(x) | \leq E_{Q_N} = |\alpha|\frac{\Gamma(N+1/2)}{2\pi\sqrt{2}|x|^{1/2+N}}\exp\left(2\alpha\rho\right)
|\mathrm{G}(x) - \mathrm{Q}_N(x) | \leq E_{\mathrm{Q}_N}(x) = \frac{(2N-1)!!}{\pi^{N+1}x^{2N+1}}.
\end{equation}
%Where $\Gamma(x)$ is the gamma function and $\alpha$, $\rho$ are functions of $x$ and given as follows
%\begin{equation*}
%\alpha(x) = \frac{2|x|}{2|x| - 1}, \qquad \rho(x) = \ds\frac{1}{4|x|} + \frac{|8x + 1|}{4|x||2x-1|^2}.
%\end{equation*}
The expansion in eq.~\eqref{eq:QN} is connected to the generalized hypergeometric series $_2F_0(1/2,1,-2\I/(\pi x^2))$. The bound in eq.~\eqref{eq:EQN} is a much sharper then directly estimating the remainder term of the hypergeometric series, see 13.7.5 from NIST \cite{NIST:DLMF}.
\section{Technique}\label{sec:4}
We start this section by making some observations. The Taylor expansion $\mathrm{T}_N$ is accurate near zero.  As $|x|$ becomes increasingly large, more terms are needed to maintain accuracy.  On the other hand, the asymptotic expansion $\mathrm{Q}_N$ is accurate for $|x|$ large and as $|x|$ become small it needs more terms to maintain accuracy.  The modified trapezoid rule does have a global bound; see \cite{alazah2014computing}, but in general is more prone to round-off errors for $|x|$ too small or too large. Given the parity of the Fresnel integrals, we consider positive $x$ and split the non-negative reals into 3 sub-intervals with cut-off points $x_1$ and $x_2$.  Orders of approximation $N_1,N_2, N_3$ are also selected for the Taylor polynomial, $\mathrm{T}_N$, the modified trapezoid rule, $\mathrm{G}_N$, and the asymptotic expansion, $\mathrm{Q}_N$, respectively. We give a methodical approach to choosing the parameters $N_1,N_2,N_3$, $x_1$, $x_2$ and we select these parameters for double precision. In doing so we balance the float point operations between the 3 approximations and maintain the accuracy of the approximation.  This provides and accurate approximation, near the desired machine precision with no increase in computational time for any $x > 0$. We remark that if one would like another precision, the parameters can be refined in the same manner as done below. 

We start with equation~\eqref{eq:EGN}, which is a global bound for the modified trapezoid rule. Choosing the parameter $N_2 = 12$ gives $E_{\mathrm{G}_{12}} \approx 1.0733e$-$17$ which is below double precision. Choosing the value $N_2 = 11$ gives $E_{\mathrm{G}_{11}} \approx 2.301e$-$16$ which is slightly above double precision, for this example we will choose $N_2 = 12$. As per the above discussion, for a given $N$, $\mathrm{T}_N$ proves to be a better approximation for small $x$, while $\mathrm{Q}_N$ is a better approximation for large $x$. We thus choose $N_1$ and $N_3$ based on the computational time. Doing so gives $N_1 = 14$, $N_3 = 12$.  In deciding the orders $N_1, N_2, N_3$, we can now choose $x_1$ and $x_2$ based off the precision we want. Considering eq.~\eqref{eq:ETN} at $x_1 = 0.688$ gives $E_{\mathrm{T}_{14}}(0.688) \approx 2.078e$-$16$. Also, considering eq.~\eqref{eq:EQN} at $x_2 =6.725$ gives $E_{\mathrm{Q}_{12}}(6.725) \approx 2.212e$-$16$. 

In summary, we have the approximation
\begin{equation}
G(x) \approx \tilde{G}(x) =  \begin{cases}
T_{14}(x), & x\in [0,0.688] \\
G_{12}(x), & x\in (0.688, 6.725)  .\\
Q_{12}(x), & x\in [6.725,\infty)
\end{cases}
\end{equation}
Using the the parity of $\mathrm{G}(x)$ the estimate $|\mathrm{G}(x) - \tilde{\mathrm{G}}(x) | < 2^{-52}$ holds for all $x\in\mR$.  So we achieve a double precision approximation that has similar computational time for all $x\in\mR$ while avoiding machine rounding errors.
%---NUMERICAL---------------------------------------------------------------------------------
\section{Numerical Experiments}\label{sec:5}
Alazah et al. (see \cite{alazah2014computing}) test the validity of their method for arguments between 0 and 1000 and endorse the truncated Taylor series for small arguments, a well known approach. Many other methods present in literature use polynomial or rational approximations, splitting the domain into several regions. As far as we know, the most accurate one is Cody \cite{cody1968chebyshev} with relative errors $\leq 10^{-15.58}$. However the form for the intermediate region and cut-off points are found experimentally. Mielenz \cite{mielenz1997computation} uses an asymptotic expansion for large arguments equivalent to a particular case of the form we use fixing $N=6$. We confirm the precision of our technique here and test its efficiency. 

Numerical experiments were run with {\scshape MatLab} ver. R2020a using a 5th generation Intel i7-5500U CPU @ 2.4GHz, with 16GB of PC3L-12800 DDR3L ram. We compare our approximation against two sources.  First the native Frensel cosine and sine functions for {\scshape MatLab}, we denote this by $\mathrm{G}_m(x)$ = {\itshape fresnelc}$(x) + \I ${\itshape fresnels}$(x)$ with the {\itshape fresnelc} and {\itshape fresnels} cosine and sine functions. We compare with an approximation based off Mielenz~\cite{mielenz1997computation}, which uses a spline interpolant. This will be denoted as $\mathrm{G}_s(x)$ and can be found in the {\scshape MatLab} file exchange\footnote{John D'Errico (2020). FresnelS and FresnelC \\
(https://www.mathworks.com/matlabcentral/fileexchange/28765-fresnels-and-fresnelc),\\
MATLAB Central File Exchange.}.

\begin{figure}[htpb!]
\subfigure[Error $|\mathrm{G}_m - \tilde{\mathrm{G}}|$]
{\begin{minipage}[t]{0.45\linewidth}
\includegraphics*[scale = .35]{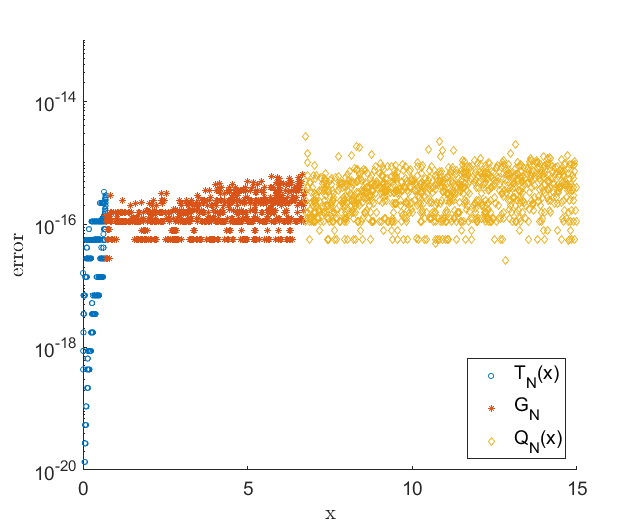}
\end{minipage}}
\subfigure[Error $|\mathrm{G}_s - \tilde{\mathrm{G}}|$]
{\begin{minipage}[t]{0.45\linewidth}
\centering
\includegraphics*[scale = .35]{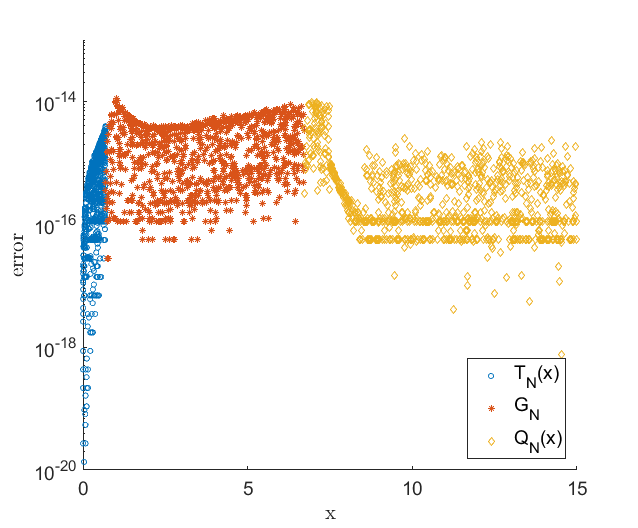}
\end{minipage}}
%\vspace{-10pt}
\caption{Error plots for $x\in[0,15]$. 1000 random test points are taken for each subinterval partitioned by $x_1 = 0.688$, $x_2 = 6.725$. (a) The native {\scshape MatLab} function $\mathrm{G}_m(x)$ = {\itshape fresnelc}$(x) + \I ${\itshape fresnels}$(x)$ vs $\tilde{\mathrm{G}}(x)$,  (b) Error plot of $\mathrm{G}_s(x)$ vs $\tilde{\mathrm{G}}(x)$.}\label{fig:1}
\end{figure}

Figure~\ref{fig:1} shows the errors between the approximation $\tilde{\mathrm{G}}(x)$, $\mathrm{G}_m(x)$ and $\mathrm{G}_s(x)$. For each subinterval, partitioned by $x_1 = 0.688$ and $x_2 = 6.725$, each approximation is evaluated at 1000 random points. The average of the errors is in Table~\ref{tab:table2}. Although $\mathrm{G}_m(x)$ and $\mathrm{G}_s(x)$ are approximations themselves, there is consistent behavior between these approximations. The error reflects the accuracy of the approximation for a the given subinterval.
\begin{table}[h!]
\begin{center}
	\begin{tabular}{c c c c}
	 &\multicolumn{3}{c}{Approximation} \\
	\cline{2-4} \vspace{-10pt} \\
		 Interval & \textbf{$\tilde{\mathrm{G}}$} & $\mathrm{G}_s$ & $\mathrm{G}_m$  \\
			\hline
		$[0,x_1]$  &3.2951313293e-05 & 1.2483252225e-04 & 2.776233718e-02  \\
		$(x_1,x_2)$&2.9528400045e-05 & 8.6929393952e-05 & 2.874295949e-02  \\
		$[x_2,15]$ &3.2216610024e-05 & 4.2383601564e-05 & 2.809127640e-02 \\
		\\
		\end{tabular}
					\caption{Average time comparison for Figure~\ref{fig:1}. 1000 random points are taken for each subinterval partitioned by $x_1 = 0.688$, $x_2 = 6.725$. The time for computing the approximation is record with commands {\itshape tic}/{\itshape toc}.}\label{tab:table1}
	\end{center}
\end{table}
\vspace{-20pt}
\begin{table}[h!]
\begin{center}
	\begin{tabular}{c c c}
	 &\multicolumn{2}{c}{Average point-wise error} \\
	\cline{2-3} \vspace{-10pt} \\
		 Interval & $|\mathrm{G}_m-\tilde{\mathrm{G}}|$ & $|\mathrm{G}_s-\tilde{\mathrm{G}}|$  \\
			\hline
		$[0,x_1]$  &4.3744537030e-17 & 6.1597220741e-16\\
		$(x_1,x_2)$&1.4010574651e-16 & 2.1708872956e-15  \\
		$[x_2,15]$ &4.0216009902e-16 & 6.8351638744e-16  \\
		\\
		\end{tabular}
					\caption{Average point-wise error comparison for Figure~\ref{fig:1}. 1000 random points are taken for each subinterval partitioned by $x_1 = 0.688$, $x_2 = 6.725$.}\label{tab:table2}
	\end{center}
\end{table}

\vspace{-10pt}

\begin{figure}[htpb!]
\subfigure[$|\mathrm{G}_m - \mathrm{Q}_N|$]
{\begin{minipage}[t]{0.45\linewidth}\label{fig:2a}
\includegraphics*[scale = .33]{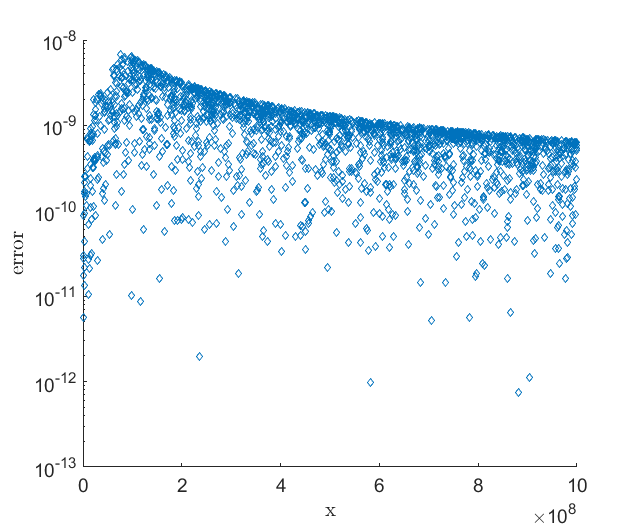}
\end{minipage}}
\subfigure[$|\mathrm{G}_s - \mathrm{G}_m|$]
{\begin{minipage}[t]{0.45\linewidth}\label{fig:2b}
\centering
\includegraphics*[scale = .33]{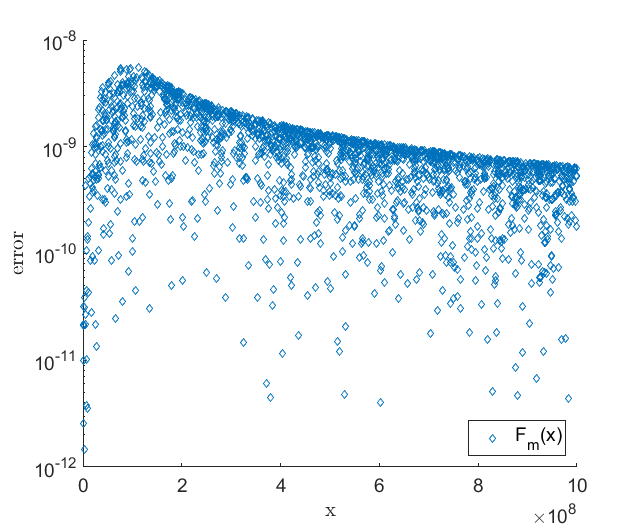}
\end{minipage}}
%\\
%\subfigure[$|\mathrm{G}_s - \mathrm{Q}_N|$]
%{\begin{minipage}[t]{0.45\linewidth}\label{fig:2c}
%\centering
%\includegraphics*[scale = .33]{image/p4_3.png}
%\end{minipage}}
%\vspace{-10pt}
\caption{Error plots for $x\in[10,1.e9]$ with 2000 random test points. (a) The native {\scshape MatLab} functions $\mathrm{G}_m(x)$ and The asymptotic expansion $\mathrm{Q}_{12}(x)$. (b) The native {\scshape MatLab} functions $\mathrm{G}_m(x)$ and $\mathrm{G}_s(x)$.}\label{fig:2} %(c) $\mathrm{G}_s(x)$ and $\mathrm{Q}_{12}(x)$, the average error is computed to be 6.2355314326e-16}\label{fig:2}
\end{figure}

 Figure~\ref{fig:2} shows a comparison for large $x$, we use 2000 random test points between the interval $[10,1.e9]$. We show a comparison between the errors generated by the native functions in {\scshape MatLab}, $\mathrm{G}_m(x)$ and the asymptotic expansion $\mathrm{Q}_{12}(x)$, the errors generated by $\mathrm{G}_m(x)$ and $\mathrm{G}_s(x)$.  This illustrates the limit of precision of the approximation $\mathrm{G}_m(x)$. We note that $\mathrm{G}_s(x)$ and $\tilde{\mathrm{G}}(x)$ use the same asymptotic expansion for large $x$, but their implementations are different.  It should also be noted that both approximations $\mathrm{G}_s(x)$ and $\tilde{\mathrm{G}}(x)$ are fast to compute, the time can be seen in Table~\eqref{tab:table1}. If one were to need a quicker implementation, one could precompute the coefficients for given parameters $N_1, N_2, N_3$.

%---CONCLUSION--------------------------------------------------------------------------------
\section{Conclusion}\label{sec:6}
In this work we have presented a systematic way to compute Fresnel integrals up to an arbitrary machine precision. Our numerical experiments verify this for double precision.  If more accuracy is required, one can select the parameters in a similar fashion to what is done in Section~\ref{sec:4}. In Section~\ref{sec:5} we are also able to justify the results from previous works where justification was heuristic. With a proper implementation, the method presented can improve the performance for computing Fresnel integrals in standard packages such as {\itshape SciPy}~\cite{2020SciPy-NMeth}. In computer aided geometric design, this technique can be used to plot clothoid segments for transition curves and clothoidal splines.

The underlying philosophy combining Taylor series, some clever form of quadrature and asymptotic expansion can be investigated for the computation of other special functions and more exotic curves arising in CAGD.  
%\section*{Acknowledgment}
%---APPENDIX----------------------------------------------------------------------------------
\appendix
%\label{sec:append}
%\nocite{*}
%---REFERENCES--------------------------------------------------------------------------------
\bibliographystyle{abbrv}
\bibliography{bibliography}
\end{document}